\documentclass{amsart}
\usepackage{amssymb}
\usepackage{amsfonts}
\usepackage{amsmath}
\usepackage{graphicx}
\usepackage[pagebackref]{hyperref}
\theoremstyle{plain}

\theoremstyle{definition}

\theoremstyle{remark}

\numberwithin{equation}{section}

\def\1{{\rm (1)}}
\def\2{{\rm (2)}}
\def\3{{\rm (3)}}
\def\4{{\rm (4)}}
\def\5{{\rm (5)}}

\begin{document}

\title[Corrigendium:  Krull dimension, Overrings and Semistar Operations of an Integral Domain]
{Corrigendium: Krull dimension, Overrings and Semistar Operations of an Integral Domain.
[J. Algebra, 321 (2009), 1497--1509]}


\author{A. Mimouni}

\address{Department of Mathematics and Statistics,
King Fahd University of Petroleum \& Minerals, P. O. Box 278,
Dhahran 31261, Saudi Arabia}

\email{amimouni@kfupm.edu.sa}

\thanks{This work is supported by KFUPM}



\maketitle
\noindent{\bf Theorem 2.9} Let $R$ be an integrally closed  domain. Then
the following conditions are equivalent.\\
$(i)$ $|O(R)|=7+dimR$.\\
$(ii)$ $|SSFc(R)|=7 + dimR$.\\
$(iii)$ $R$ is a Pr\"ufer domain, and either\\
$(a)$ $R$ has exactly three maximal ideals $M_{1}, M_{2}, M_{3}$ and $Spec(R)=\{0=P_{0}\subsetneq P_{1}\subsetneq \dots\subsetneq
P_{r-1}\subsetneq M_{1}\cap M_{2}\cap M_{3}\}$; or\\
$(b)$ $dimR\geq 5$, $R$ has exactly two maximal ideals $M$ and $N$
and $Spec(R)=\{0=P_{0}\subsetneq P_{1}\subsetneq \dots\subsetneq
P_{r-1}\subsetneq P_{r}=M, N, P_{r-4}\nsubseteq N,
P_{r-5}\subsetneq N\}$; or\\
$(c)$ $dimR\geq 2$, $R$ has exactly two maximal ideals $M$ and $N$
and $Spec(R)=\{0=P_{0}\subsetneq P_{1}\subsetneq \dots\subsetneq
P_{r-2}\subsetneq P_{r-1}\subsetneq P_{r}=M, P_{r-2}\subsetneq
Q\subsetneq N\}$.

\begin{proof} Only Case 1 need a slight modification proof.\\
\noindent{\bf Case 1} $|Max(R)|=3$. Assume that $Max(R)=\{M_{1}, M_{2}, M_{3}\}$ and set
$htM_{1}=dimR=r\geq 2$. Let $0=P_{0}\subsetneq P\subsetneq M_{1}$ be
a chain of prime ideals of $R$. We claim that $P_{r-1}\subseteq M_{2}\cap M_{3}$. Indeed, 
suppose that $P_{r-1}\nsubseteqq M_{2}$. Then  $O(R)\supseteq \{R, L, R_{P_{1}},\dots, R_{P_{r-1}}, R_{M_{1}}, R_{M_{2}}, R_{M_{3}}, R_{M_{1}}\cap R_{M_{2}}, R_{M_{1}}\cap R_{M_{3}}, R_{M_{2}}\cap R_{M_{3}}, R_{P_{r-1}}\cap R_{M_{2}}\}$. Hence
$7+r=7+dimR=|O(R)|\geq 9+r=9+dimR$, which is absurd. Thus the assertions $(a)$ is satisfied and $spec(R)$ is of the form.\\
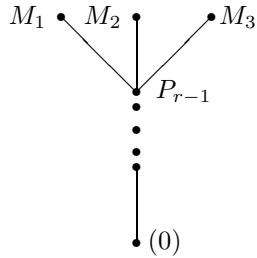
\[\setlength{\unitlength}{1mm}
\begin{picture}(80,20)(-10,00)
\put(40,20){\line(0,-1){10}} \put(40,10){\line(1,1){10}}
\put(40,10){\line(-1,1){10}}
\put(40,00){\line(0,-1){10}} 
\put(40,08){\circle*{1}}
\put(40,05){\circle*{1}}
\put(40,02){\circle*{1}}
\put(40,-10){\circle*{1}}
\put(40,20){\circle*{1}} \put(38,20){\makebox(0,0)[r]{$M_{2}$}}
\put(30,20){\circle*{1}} \put(28,20){\makebox(0,0)[r]{$M_{1}$}}
\put(50,20){\circle*{1}} \put(56,20){\makebox(0,0)[r]{$M_{3}$}}
\put(40,10){\circle*{1}} \put(50,10){\makebox(0,0)[r]{$P_{r-1}$}}
\put(40,00){\circle*{1}}\put(46,-10){\makebox(0,0)[r]{$(0)$}}
\end{picture}\]
\end{proof}
\bigskip

\noindent{\bf Example 3.4.} A one-dimensional Noetherian local domain $R$
 such $|O(R)|=3+dimR$ and $|SSFc(R)|=5+dimR$.\\
Let $k$ be a field and $X$ an indeterminate over $k$. Set
$R=k[[X^{3}, X^{4}, X^{5}]]$. Then $R$ is a one-dimensional
Noetherian local domain which is not divisorial (Lemma 3.1,
since $M^{-1}=k[[X]]$ and $R\subsetneq k[[X^{2}, X^{3}]]\subsetneq
M^{-1}$). Also it is easy to see that $O(R)=\{R\subsetneq
R_{1}=k[[X^{2}, X^{3}]]\subsetneq R_{2}=R'=k[[X]]\subsetneq
L=qf(R)=k((X))\}$ and each proper overring of $R$ is divisorial
(Lemma 3.1) so $fgv$ domain. Clearly $|O(R)|=4=3+dimR$,
however $|SSFc(R)|=6=5+dimR$ (in fact $SSFc(R)=\{e, \bar{d}, \overline{*_{R_{1}}},
\bar{t}, \{*_{R_{i}}\}_{\{1\leq i\leq 2\}}\}$, where $E^{\overline{*_{R_{1}}}}=ER_{1}\cap E_{v}$).\\

\medskip
\noindent{\bf Example 3.5.} A one-dimensional Noetherian local
domain with $|O(R)|=|SSFc(R)|\\=5+dimR$.\\
Let $k$ be a field and $X$ an indeterminate over $k$. Set
$R=k[[X^{2}, X^{9}]]$. Then $R$ is a one-dimensional Noetherian
local domain which is divisorial (since $M^{-1}=k[[X^{2}, X^{7}]])$. Also it is easy to see that $O(R)=\{R\subsetneq
R_{1}=k[[X^{2}, X^{7}]]\subsetneq R_{2}=k[[X^{2}, X^{5}]]\subsetneq R_{3}=k[[X^{2}, X^{3}]]\subsetneq
R_{4}=R'=k[[X]]\subsetneq L=qf(R)=k((X))\}$ and each overring of $R$
is divisorial (Lemma 3.1) so $fgv$ domain. By
Theorem 2.3, $|SSFc(R)|=|O(R)|=6=5+dimR$ (in fact
$SSFc(R)=\{e, \bar{d},\\ \{*_{R_{i}}\}_{\{1\leq i\leq 4\}}\}$).\\

\medskip
\noindent{\bf Example 3.6.}\label{ECE.6} The following is an example of a domain
$R$ with $|O(R)|=5+dimR$ and $|SSFc(R)|=\infty$.\\
Let $\mathbb{Q}$ be the field of rational numbers, $X$ an
indeterminate over $\mathbb{Q}$. Set $V=\mathbb{Q}(\sqrt{2},
\sqrt{3})[[X]]= \mathbb{Q}(\sqrt{2}, \sqrt{3})+M$, where $M=XV$ and
$R=\mathbb{Q}+M$. By [2, Theorem 2.1], $R$ is a
one-dimensional Noetherian local domain with maximal ideal $M$.
Since each overring of $R$ is comparable to $V$ ([2, Theorem
3.1]), then it is easy to see that $O(R)=\{R, V, L=qf(R),
R_{1}=\mathbb{Q}(\sqrt{2})+M, R_{2}=\mathbb{Q}(\sqrt{3})+M,
R_{3}=\mathbb{Q}(\sqrt{6})+M\}$ and all proper overrings of $R$ are
divisorial ([2, Corollary 4.4]). Hence $|O(R)|=6=5+dimR$. Now, since $[\mathbb{Q}(\sqrt{2},
\sqrt{3}):\mathbb{Q}]=4$ and $Q$ is infinite, $|S(R)|=\infty$ and therefore $|SSFc(R)|=\infty$.\\

\medskip
\noindent{\bf Example 3.7.} A one-dimensional Noetherian local domain $R$ such that
$|SSFc(R)|=|O(R)|=6+dimR$.\\
Let $k$ be a field and $X$ an indeterminate over $k$. Let
$R=k[[X^{2}, X^{11}]]$. Since $M^{-1}=k[[X^{2}, X^{9}]]$, by
Lemma 3.1, $R$ is divisorial. Now, it is easy to see that
$O(R)=\{R\subsetneq R_{1}=k[[X^{2}, X^{9}]]\subseteq R_{2}=k[[X^{2}, X^{7}]]\subsetneq R_{3}=k[[X^{2}, X^{5}]]\subsetneq R_{4}=k[[X^{2}, X^{3}]]\subsetneq R_{5}=R'=k[[X]]\subsetneq L=qf(R)=k((X))\}$. Also it is easy to check that
each overring is divisorial (Lemma 3.1). By
Theorem 2.3, $|SSFc(R)|=|O(R)|=7=6+dimR$.\\

\medskip
\noindent{\bf Example 3.9.} Just change $R_{1}=k[[X^{3}, X^{5}]]$ to $R_{1}=k[[X^{2}, X^{9}]]$.

\end{document}